\documentclass{article}

\usepackage[latin1]{inputenc}
\usepackage{amsmath,amssymb,amsthm}
\usepackage{graphicx,graphics,color,epsfig}
\usepackage{hyperref}

\newcommand{\finprv}{\scriptsize ~$\blacksquare$ \normalsize}

\newcommand{\veps}{\varepsilon}

\newcommand{\mpdiv}{\text{$/$\hspace{-1.42mm}\tiny${}^\circ $\normalsize}}

\newcommand{\CG}{\mathcal G}

\theoremstyle{plain}
\newtheorem{theo}{Theorem}

\newtheorem{coro}{Corollary}
\newtheorem{lemm}{Lemma}
\theoremstyle{definition}

\theoremstyle{remark}
\newtheorem{rema}{Remark}

\title{Solving the additive eigenvalue problem associated to a dynamics
   of a 2D-traffic system}
\author{Nadir Farhi\thanks{Current address: University of Texas at
Dallas, 800 West Campbell Road, Richardson, TX 75080, USA.
\texttt{nadir.farhi@utdallas.edu}}\\
\small{INRIA - Paris - Rocqencourt}\\
Domaine de Voluceau, 78153, Le Chesnay, Cedex France.}
\date{}

\begin{document}
\maketitle

\begin{abstract}
  This is a technical note where we solve the
  additive eigenvalue problem
  associated to a dynamics of a 2D-traffic
  system. The traffic modeling is not explained here. It is
  available in~\cite{Far08}. It consists of a microscopic road traffic
  model of two circular roads crossing on one junction managed with the
  priority-to-the-right rule. It is based on Petri nets and minplus algebra.
  One of our objectives in~\cite{Far08} was to derive the fundamental diagram of 2D-traffic,
  which is the relation between the density and the flow of
  vehicles. The dynamics of this system, derived from a Petri net
  design, is non monotone and additively homogeneous of degree 1.
  In this note, we solve the additive eigenvalue problem associated to this dynamics.
\end{abstract}

%--------------------------
\section{Introduction}
%--------------------------

In this note we solve the additive eigenvalue problem (or the
time-independent system) associated to the dynamics of a basic
2D-traffic model considered in~\cite{Far08}. It is a system of two
circular roads crossing on one junction managed with the
\emph{priority-to-the-right} rule. The model is based on Petri
nets and on minplus algebra~\cite{BCOQ92}. It is an extension to
an existing 1D-traffic model~\cite{FGQ05,LMQ01}, which gives the
average speed of vehicles on one circular road as an eigenvalue of
a minplus matrix, and thus allows the derivation of the
fundamental diagram of 1D-traffic (the relation between the
density and the flow of vehicles on the road).

We give a solution of the eigenvalue problem. We show that the
eigenvalue is not necessarily unique, but is given in terms of two
main quantities which are interpreted in terms of traffic as the
density $d$ of vehicles in the system, and a parameter $r$ giving
the ratio of the non priority road size with respect to the size
of the whole system. We give a condition on the parameter $r$ such
that the eigenvalue is unique and positive for non-high densities.
In this case the eigenvalue problem can be written as a dynamic
programming equation of a stochastic optimal control problem.

We use the minplus algebra notations, mainly for reason of
compactness but also to use some classical results of this
algebra~\cite{BCOQ92}. In addition, the following notations are
also used: $a\mpdiv b$ denotes $a-b$, $\sqrt{a}$ denotes $a/2$,
and $b^a$ denotes $ab$.

%The notations we use below are summarized in~:
%\begin{itemize}
%  \item For $A$ and $B$ two square matrices, $A\oplus B$
%    denotes the elementwise minimum of $A$ and~$B$.
%  \item For a $(n\times m)$ matrix $A$ and a $(m\times p)$ matrix $B$,
%    the $(n\times p)$ matrix $A\otimes B$ is defined~:
%    $$(A\otimes B)_{ij}=\min_{1\leq k\leq m}(A_{ik}+B_{kj}).$$
%  \item For $a$ and $b$ two real numbers, $a\mpdiv b$ denotes
%    $a-b$, and $\sqrt{a}$ denotes $a/2$.
%  \item $\veps$ and $e$ used as scalars denote $+\infty$ and $0$
%    respectively. Used as matrices, $\veps$ denotes the matrix whose
%    components are all equal to $+\infty$, and $e$ the square
%    matrix whose components are all equal to $+\infty$ except the
%    diagonal ones which are equal to $0$.
%\end{itemize}
The traffic dynamics is the following (see~\cite{Far08}):
\begin{align}
  &x_q^{k+1}=a_{q-1}x_{q-1}^k\oplus\bar{a}_qx_{q+1}^k,\; \label{eq}\\
  &\quad \quad \quad q\in\{2,\ldots,n-1,n+2,\ldots,n+m-1\} \;, \nonumber\\
  &x_{n}^{k+1}=\bar{a}_n x_1^k x_{n+1}^k\mpdiv x_{n+m}^{k+1}\oplus a_{n-1}x_{n-1}^k\;,\label{nimplicit}\\
  (DS):\quad \quad \quad &x_{n+m}^{k+1}=\bar{a}_{n+m} x_1^k x_{n+1}^k\mpdiv x_{n}^k\oplus a_{n+m-1} x_{n+m-1}^k\;,
    \label{implicit}\\
  &x_{1}^{k+1}=a_{n+m}\sqrt{x_{n}^kx_{n+m}^k}\oplus\bar{a}_1 x_2^k\;,\label{ar1}\\
  &x_{n+1}^{k+1}=a_{n}\sqrt{x_{n}^k x_{n+m}^k}\oplus \bar{a}_{n+1}
    x_{n+2}^k\;,\label{ar2}
\end{align}
with the (traffic) constraints:
\begin{equation}\label{cont}
  \begin{cases}
    0\leq a_i\leq 1 & i=1,2,\ldots, n+m,\\
    \bar{a}_i=1\mpdiv a_i & i\neq n, n+m,\\
    \bar{a}_n=\bar{a}_{n+m}=1\mpdiv (a_n a_{n+m}),\\
    a_n a_{n+m}\leq 1.
  \end{cases}
\end{equation}

For example, in the usual algebra, the equation (\ref{nimplicit})
is written:
$$x_{n}^{k+1}=\min\left\{\bar{a}_n+x_1^k+x_{n+1}^k-x_{n+m}^{k+1}\;,\;a_{n-1}+x_{n-1}^k\right\},$$
when the equation (\ref{ar1}) is written:
$$x_1^{k+1}=\min\left\{a_{n+m}+\frac{x_{n}^k+x_{n+m}^k}{2}\;,\;\bar{a}_1 +x_2^k\right\}.$$
This system of equations is implicit but it is triangular, so its
trajectory is unique.

We denote by $d$ the following quantity (which is interpreted in
terms of traffic as the density of vehicles in the system):
\begin{equation}\label{nota}
  d=\frac{1}{n+m-1}\sum_{i=1}^{n+m} a_i.
\end{equation}

%-----------------------------------------------------
\section{Solving the additive eigenvalue problem}
%-----------------------------------------------------

The additive eigenvalue problem corresponding to the dynamics (DS)
is:
\begin{align}
  &\lambda x_i=a_{i-1}x_{i-1}\oplus\bar{a}_ix_{i+1},\; \label{EV-eq}\\
  &\quad \quad \quad i\in\{2,\ldots,n-1,n+2,\ldots,n+m-1\} \;, \nonumber\\
  &\lambda x_{n}=\bar{a}_n x_1 x_{n+1}\mpdiv (\lambda x_{n+m})\oplus a_{n-1}x_{n-1}\;,\label{EV-nimplicit}\\
  (EV):\quad \quad \quad
  &\lambda x_{n+m}=\bar{a}_{n+m} x_1 x_{n+1}\mpdiv x_n\oplus a_{n+m-1} x_{n+m-1}\;,
    \label{EV-implicit}\\
  &\lambda x_1=a_n\sqrt{x_n x_{n+m}}\oplus\bar{a}_1 x_2\;,\label{EV-ar1}\\
  &\lambda x_{n+1}=a_{n+m}\sqrt{x_n x_{n+m}}\oplus \bar{a}_{n+1} x_{n+2}\;,\label{EV-ar2}
\end{align}
with the constraints (\ref{cont}) and the notation (\ref{nota}).
Our aim in this note is to solve the system $(EV)$.
\begin{theo}
  Solving the system $(EV)$ is equivalent to solving the following simplified system (SS):
  \begin{align}
    &x_i=(a_{i-1}\mpdiv\lambda) x_{i-1}\oplus (\bar{a}_i\mpdiv\lambda)x_{i+1},\; \label{SS-eq}\\
    &\quad \quad \quad i\in\{2,\ldots,n-1,n+2,\ldots,n+m-1\} \;, \nonumber\\
    &x_n=(\bar{a}_n\mpdiv\lambda^2) x_1 x_{n+1}\mpdiv x_{n+m} \oplus (b_n\mpdiv\lambda^{n-1}) x_1\;,
       \label{SS-nimplicit}\\
    (SS):\quad \quad \quad
    &x_{n+m}=(\bar{a}_{n+m}\mpdiv\lambda) x_1 x_{n+1}\mpdiv x_n \oplus (b_m\mpdiv\lambda^{m-1}) x_{n+1}\;,
       \label{SS-implicit}\\
    &x_1=(a_n\mpdiv\lambda) \sqrt{x_n x_{n+m}} \oplus (\bar{b}_n\mpdiv\lambda^{n-1}) x_n\;,\label{SS-ar1}\\
    &x_{n+1}=(a_{n+m}\mpdiv\lambda) \sqrt{x_n x_{n+m}} \oplus (\bar{b}_m\mpdiv\lambda^{m-1}) x_{n+m}\;,
       \label{SS-ar2}
  \end{align}
  where $b_n=\bigotimes_{i=1}^{n-1}a_i$, $\bar{b}_n=\bigotimes_{i=1}^{n-1}\bar{a}_i$,
        $b_m=\bigotimes_{i=n+1}^{n+m-1}a_i$ and $\bar{b}_m=\bigotimes_{i=n+1}^{n+m-1}\bar{a}_i$.
\end{theo}
\proof We proceed in two steps:
\begin{itemize}
  \item First we show that if $(\lambda,x)$ is a solution of the system $(EV)$ then $\lambda\leq 1/4$.
  Indeed, from the equations (\ref{EV-nimplicit}), (\ref{EV-ar1}) and (\ref{EV-ar2}), we have:
  \begin{align}
    &\lambda x_{n}=\bar{a}_n x_1 x_{n+1}\mpdiv (\lambda x_{n+m})\oplus a_{n-1}x_{n-1}, \nonumber \\
    &\lambda x_1\leq a_n\sqrt{x_n x_{n+m}}, \nonumber \\
    &\lambda x_{n+1}\leq a_{n+m}\sqrt{x_n x_{n+m}}. \nonumber
  \end{align}
  Then by multiplying (standard adding) the terms of these inequalities, we obtain $\lambda^4\leq 1$, since
  $\bar{a}_{n+m}a_n a_{n+m}=1$.

  \item We see that if $n=m=2$ the systems $(EV)$ and $(SS)$ coincide.
    For $n$ and $m$ fixed, we denote by $EV(n,m)$ and $SS(n,m)$ the corresponding systems.
    By induction on $n$ and $m$, we suppose that $EV(n,m)\Leftrightarrow SS(n,m)$, and we show
    that $EV(n+1,m)\Leftrightarrow SS(n+1,m)$ and $EV(n,m+1)\Leftrightarrow SS(n,m+1)$.
    \begin{itemize}
      \item To show that $EV(n+1,m)\Leftrightarrow SS(n+1,m)$, we
        eliminate the variable $x_n$ in $EV(n+1,m)$ which gives a
        system $EV(n,m)$, then we use the induction assumption.

        Indeed the problem $EV(n+1,m)$ is written as follows:
        \begin{align}
          &\lambda x_i=a_{i-1}x_{i-1}\oplus\bar{a}_ix_{i+1},\; \label{I-eq}\\
          &\quad \quad \quad i\in\{2,\ldots,n,n+3,\ldots,n+m\} \;, \nonumber\\
          &\lambda x_{n+1}=\bar{a}_{n+1} x_1 x_{n+2}\mpdiv (\lambda x_{n+1+m})\oplus a_n x_{n}\;,\label{I-nimplicit}\\
          %(EV):\quad \quad \quad
          &\lambda x_{n+1+m}=\bar{a}_{n+1+m} x_1 x_{n+2}\mpdiv x_{n+1}\oplus a_{n+m} x_{n+m}\;,
            \label{I-implicit}\\
          &\lambda x_1=a_{n+1}\sqrt{x_{n+1} x_{n+1+m}}\oplus\bar{a}_1 x_2\;,\label{I-ar1}\\
          &\lambda x_{n+2}=a_{n+1+m}\sqrt{x_{n+1} x_{n+1+m}}\oplus \bar{a}_{n+2} x_{n+3}\;,\label{I-ar2}
        \end{align}

        Using the expression of $x_n$ in (\ref{I-eq}), we replace it in the expression of $x_{n+1}$ in
        (\ref{I-nimplicit}). We obtain:
        $$\lambda x_{n+1}=\bar{a}_{n+1} x_1 x_{n+2}\mpdiv (\lambda x_{n+1+m})\oplus
                                a_n [(a_{n-1}\mpdiv \lambda)x_{n-1}\oplus (\bar{a}_n\mpdiv \lambda)x_{n+1}],$$
        which gives:
        \begin{equation}\label{eqnp1}
           \lambda x_{n+1}=\bar{a}_{n+1} x_1 x_{n+2}\mpdiv (\lambda x_{n+1+m})\oplus
                                   (a_n a_{n-1}\mpdiv \lambda)x_{n-1},
        \end{equation}
        because $\lambda x_{n+1}<(a_n\bar{a}_n\mpdiv
        \lambda)x_{n+1}$ since $\lambda\leq 1/4<1/2$ and
        $a_n\bar{a}_n=1$.

        Using the expression of $x_n$ in (\ref{I-eq}), we replace
        it in the expression of $x_{n-1}$ in (\ref{I-eq}) also. We  obtain:
        $$\lambda x_{n-1}=a_{n-2} x_{n-2}\oplus \bar{a}_{n-1} [(a_{n-1}\mpdiv\lambda)x_{n-1}\oplus
                                   (\bar{a}_n\mpdiv\lambda)x_{n+1}],$$
        which gives:
        \begin{equation}\label{eqnm1}
            \lambda x_{n-1}=a_{n-2} x_{n-2}\oplus (\bar{a}_{n-1}\bar{a}_n\mpdiv\lambda)x_{n+1},
        \end{equation}
        because $\lambda x_{n-1}<(\bar{a}_{n-1}a_{n-1}\mpdiv\lambda)x_{n-1}$, since $\lambda\leq
        1/4<1/2$ and $\bar{a}_{n-1}a_{n-1}=1$.

        Let us denote by $\mathcal I_1$ and $\mathcal I_2$ the following sets of indexes:
        \begin{align}
          & \mathcal I_1=\{1\leq i\leq n+1+m,\;\;i\neq 1,n+1,n+2,n+1+m\}, \nonumber \\
          & \mathcal I_2=\{1\leq i\leq n+1+m,\;\; i\neq 1,n-1,n,n+1,n+2,n+1+m\}. \nonumber
        \end{align}
        Thus we can conclude the following equivalence:
        \begin{equation}
          (\ref{I-eq})_{i\in\mathcal I_1} \Leftrightarrow \big[(\ref{I-eq})_{i\in \mathcal
          I_2},\;(\ref{eqnm1}),\;(\ref{I-eq})_{i=n}\big].
        \end{equation}

        The equations $(\ref{I-eq})_{i\in \mathcal I_2}$ combined with the equations
        (\ref{eqnm1}), (\ref{eqnp1}), (\ref{I-implicit}),
        (\ref{I-ar1}) and (\ref{I-ar2}) form the following $EV(n,m)$ system of
        variables $x_i,\;1\leq i\leq n+1+m$ and $i\neq n$:
        \begin{align}
          &\lambda x_i=a_{i-1}x_{i-1}\oplus\bar{a}_ix_{i+1},\; \label{II-eq}\\
          &\quad \quad \quad i\in\{2,\ldots,n-2,n+3,\ldots,n+m\} \;, \nonumber\\
          &\lambda x_{n-1}=a_{n-2} x_{n-2}\oplus
          (\bar{a}_{n-1}\bar{a}_n\mpdiv\lambda)x_{n+1}\;, \label{II-eqnm1}\\
          &\lambda x_{n+1}=\bar{a}_{n+1} x_1 x_{n+2}\mpdiv (\lambda x_{n+1+m})\oplus
                                   (a_n a_{n-1}\mpdiv \lambda)x_{n-1}\;,\label{II-nimplicit}\\
          %(EV):\quad \quad \quad
          &\lambda x_{n+1+m}=\bar{a}_{n+1+m} x_1 x_{n+2}\mpdiv x_{n+1}\oplus a_{n+m} x_{n+m}\;,
            \label{II-implicit}\\
          &\lambda x_1=a_{n+1}\sqrt{x_{n+1} x_{n+1+m}}\oplus\bar{a}_1 x_2\;,\label{II-ar1}\\
          &\lambda x_{n+2}=a_{n+1+m}\sqrt{x_{n+1} x_{n+1+m}}\oplus \bar{a}_{n+2} x_{n+3}\;.\label{II-ar2}
        \end{align}
        By using the induction assumption, this later system is equivalent to the following system:
        \begin{align}
          &x_i=(a_{i-1}\mpdiv\lambda) x_{i-1}\oplus (\bar{a}_i\mpdiv\lambda)x_{i+1},\; \label{SS-I-eq}\\
          &\quad \quad \quad i\in\{2,\ldots,n-1,n+3,\ldots,n+m\} \;, \nonumber\\
          &x_{n+1}=(\bar{a}_{n+1}\mpdiv\lambda^2) x_1 x_{n+2}\mpdiv x_{n+1+m} \oplus (b_{n+1}\mpdiv\lambda^{n}) x_1\;,
               \label{SS-I-nimplicit}\\
          &x_{n+1+m}=(\bar{a}_{n+1+m}\mpdiv\lambda) x_1 x_{n+2}\mpdiv x_{n+1} \oplus (b_m\mpdiv\lambda^{m-1}) x_{n+2}\;,
               \label{SS-I-implicit}\\
          &x_1=(a_{n+1}\mpdiv\lambda) \sqrt{x_{n+1} x_{n+1+m}} \oplus (\bar{b}_{n+1}\mpdiv\lambda^{n}) x_{n+1}\;,\label{SS-I-ar1}\\
          &x_{n+2}=(a_{n+1+m}\mpdiv\lambda) \sqrt{x_{n+1} x_{n+1+m}} \oplus (\bar{b}_m\mpdiv\lambda^{m-1}) x_{n+1+m}\;.
               \label{SS-I-ar2}
        \end{align}
        By adding the equation $(\ref{I-eq})_{i=n}$ to the later
        system, we obtain $SS(n+1,m)$.
      \item We show with the same manipulations that $EV(n,m+1)\Leftrightarrow SS(n,m+1)$\finprv
    \end{itemize}
\end{itemize}
\begin{lemm}\label{thmp}(Baccelli et al.~\cite{BCOQ92})
  Given $A$ a $(m\times m)$ minplus matrix, if the weights of all the
  circuits of the graph $\CG(A)$ associated to $A$ are positive, then the equation $x=A\otimes x\oplus
  b$ admits a unique solution $x=A^*\otimes b$, where
  $$A^*=\bigoplus_{n=0}^{\infty}A^n=\bigoplus_{n=0}^{m-1}A^n\;.$$
\end{lemm}
\begin{coro}
  Solving the system $(EV)$ is equivalent to solve the following
  system of four variables:
  \begin{align}
    &x_n=(\bar{a}_n\mpdiv\lambda^2) x_1 x_{n+1}\mpdiv x_{n+m} \oplus (b_n\mpdiv\lambda^{n-1}) x_1\;,
       \label{S-nimplicit}\\
    (S):\quad \quad \quad
    &x_{n+m}=(\bar{a}_{n+m}\mpdiv\lambda) x_1 x_{n+1}\mpdiv x_n \oplus (b_m\mpdiv\lambda^{m-1}) x_{n+1}\;,
       \label{S-implicit}\\
    &x_1=(a_n\mpdiv\lambda) \sqrt{x_n x_{n+m}} \oplus (\bar{b}_n\mpdiv\lambda^{n-1}) x_n\;,\label{S-ar1}\\
    &x_{n+1}=(a_{n+m}\mpdiv\lambda) \sqrt{x_n x_{n+m}} \oplus (\bar{b}_m\mpdiv\lambda^{m-1}) x_{n+m}\;,
       \label{S-ar2}
  \end{align}
\end{coro}
\proof Taking into account the equivalence $(EV)\Leftrightarrow
(SS)$, suppose that $\lambda$, $x_1$, $x_n$, $x_{n+1}$, $x_{n+m}$
are known. To determine the other variables i.e.
$x_2,\ldots,x_{n-1}$, $x_{n+1},\ldots,x_{n+m-1}$ we have to solve
the system of equation (\ref{EV-eq}). This system is written:
\begin{equation}\label{eq-aff}
  x=A\otimes x\oplus b,
\end{equation}
where
$$x={}^t(x_2,\ldots,x_{n-1},x_{n+1},\ldots,x_{n+m-1}),$$
$$A=(e\mpdiv \lambda)
    \begin{bmatrix}
       A_1 & \veps\\
       \veps & A_2
    \end{bmatrix},$$
with
$$A_1=\begin{bmatrix}
        \veps & \bar{a}_2 & \veps & \cdots & \veps\\
        a_2 & \veps & \bar{a}_3 & \cdots & \veps\\
        \veps & a_3 & \veps & \ddots & \ddots\\
        \vdots & \veps & \ddots & \ddots & \bar{a}_{n-1}\\
        \veps & \veps & \cdots & a_{n-1} & \veps
   \end{bmatrix},$$
and
$$A_2=\begin{bmatrix}
        \veps & \bar{a}_{n+1} & \veps & \cdots & \veps\\
        a_{n+1} & \veps & \bar{a}_{n+2} & \cdots & \veps\\
        \veps & a_{n+2} & \veps & \ddots & \ddots\\
        \vdots & \veps & \ddots & \ddots & \bar{a}_{n+m-1}\\
        \veps & \veps & \cdots & a_{n+m-1} & \veps
      \end{bmatrix},$$
%$$A=(e\mpdiv \lambda)
%   \begin{bmatrix}
%        \veps & \bar{a}_2 & \veps & \cdots & \veps & \veps & \veps & \veps & \veps & \veps\\
%        a_2 & \veps & \bar{a}_3 & \cdots & \veps & \veps & \veps & \veps & \veps & \veps\\
%        \veps & a_3 & \veps & \ddots & \ddots & \veps & \veps & \veps & \veps & \veps\\
%        \vdots & \veps & \ddots & \ddots & \bar{a}_{n-1} & \veps & \veps & \veps & \veps & \veps\\
%        \veps & \veps & \cdots & a_{n-1} & \veps & \veps & \veps & \veps & \veps & \veps\\
%        \veps & \veps & \veps & \veps & \veps & \veps & \bar{a}_{n+1} & \veps & \cdots & \veps\\
%        \veps & \veps & \veps & \veps & \veps & a_{n+1} & \veps & \bar{a}_{n+2} & \cdots & \veps\\
%        \veps & \veps & \veps & \veps & \veps & \veps & a_{n+2} & \veps & \ddots & \ddots\\
%        \veps & \veps & \veps & \veps & \veps & \vdots & \veps & \ddots & \ddots & \bar{a}_{n+m-1}\\
%        \veps & \veps & \veps & \veps & \veps & \veps & \veps & \cdots & a_{n+m-1} & \veps
%   \end{bmatrix},$$
and
$$b={}^t[a_1x_1\;,\;e\;,\;\ldots\;,\;e\;,\;\bar{a}_{n-1}x_n\;,\;\;
      a_{n+1}x_{n+1}\;,\;e\;,\;\ldots\;,\;e\;,\;\bar{a}_{n+m-1}x_{n+m}].$$
Using Lemma \ref{thmp}, all the circuits of the graph associated
to the matrix $A$ have the average weight $1\mpdiv \lambda^2$
which is positive because $\lambda\leq 1/4 < 1/2$. Thus the matrix
$A^*$ exists and the solution of the system (\ref{eq-aff}) is
given by: $x=A^*\otimes b$\finprv

%-------------------------------------
\section{Solving the system $(S)$:}
%-------------------------------------

\begin{figure}[htbp]
  \begin{center}
    \includegraphics[width=4.5cm,height=3.5cm]{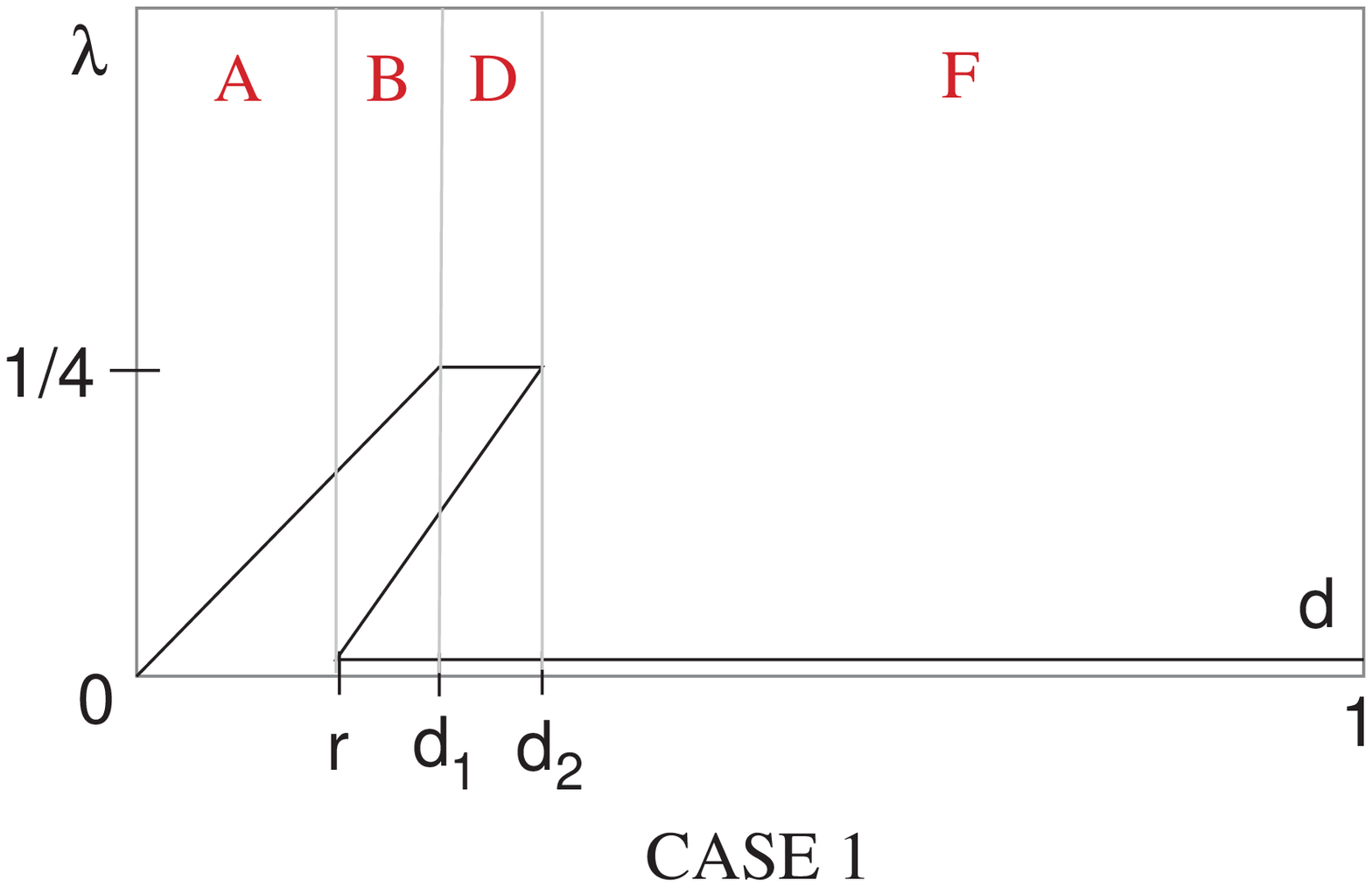}\hspace{5mm}
    \includegraphics[width=4.5cm,height=3.5cm]{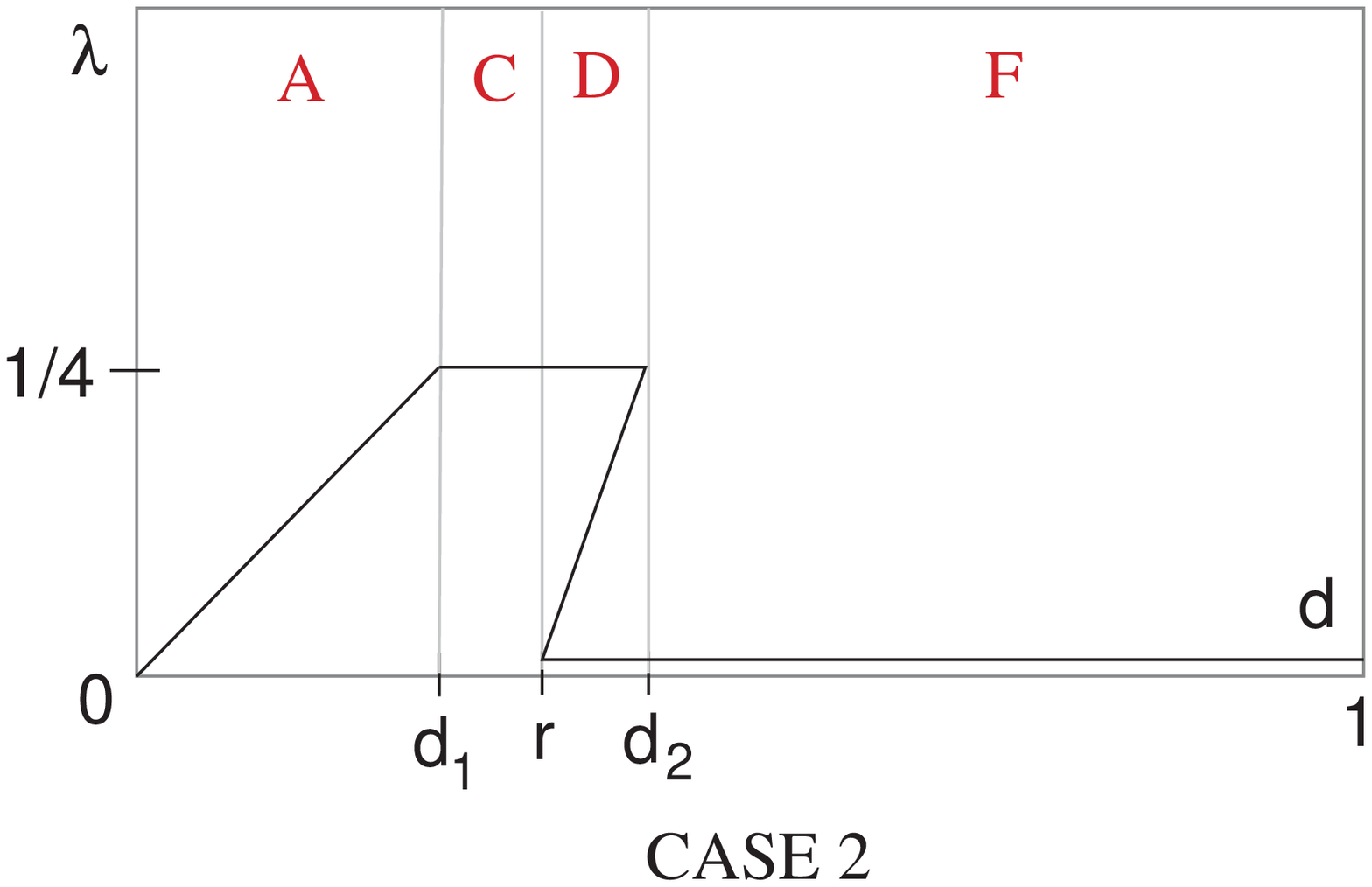}\hspace{5mm}
    \includegraphics[width=4.5cm,height=3.5cm]{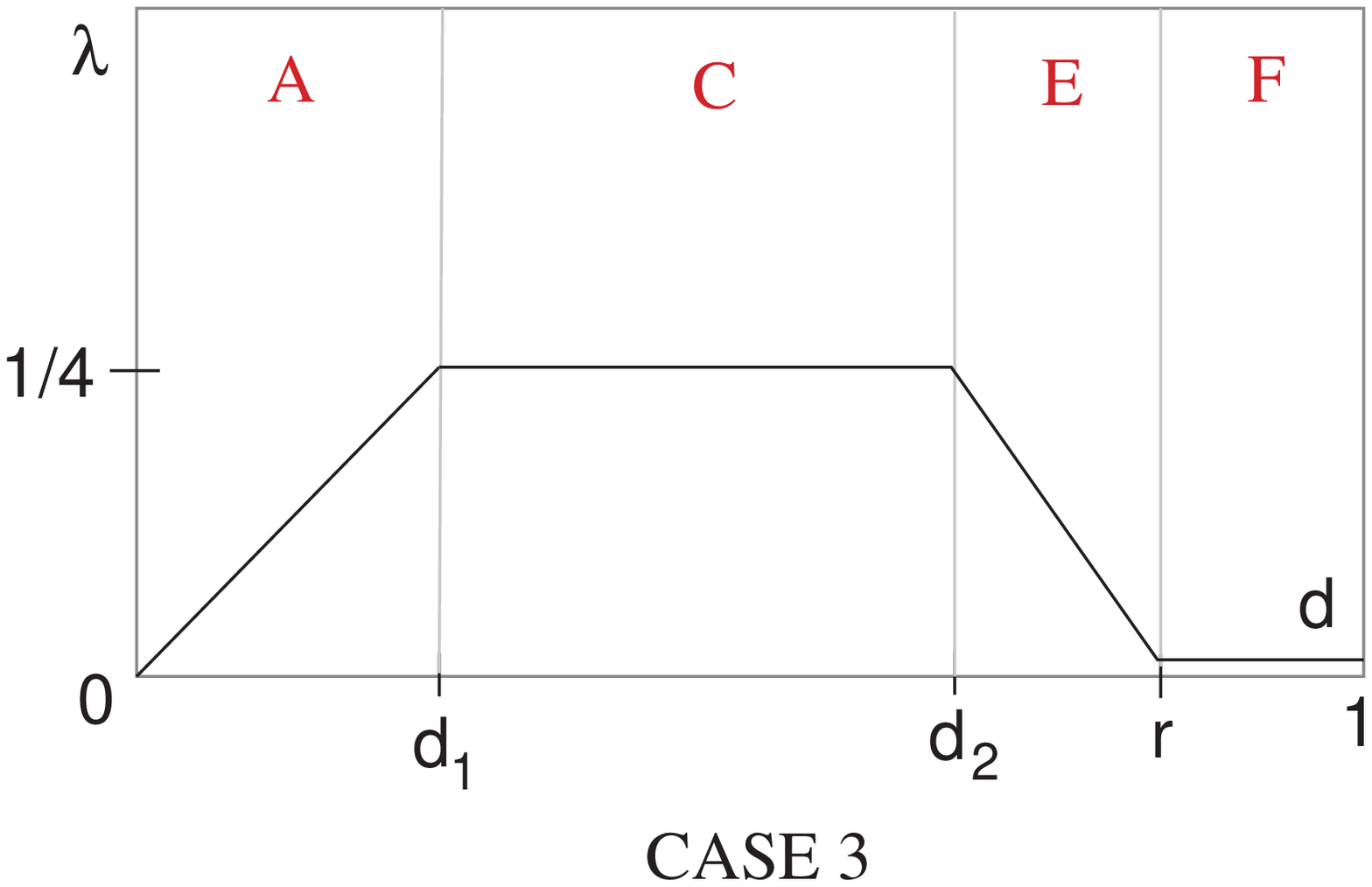}
    \caption{The curve of $\lambda$ given in Theorem~\ref{eigen} depending on $d$.}
    \label{diag}
  \end{center}
\end{figure}

Let us use the notations:
\begin{itemize}
  \item[] $r=n/(n+m-1)\;,$
  \item[] $\rho=1/(n+m-1)=r/n\;.$
\end{itemize}
\begin{theo}\label{eigen}
  There exists a solution $(\lambda,x)$ of $S$ such that $\lambda$ satisfies:
  $$0=\max\left\{\min \left\{d-(1+\rho)\lambda,\;\frac{1}{4}-\lambda,\;
      r-d-\left(2r-1+\rho\right)\lambda\right\},\;-\lambda \right\}.$$
\end{theo}
\begin{rema} Before we give the proof of Theorem~\ref{eigen},
let us explain it. Using the notations:
\begin{itemize}
  \item[] $d_1=(n+m)/[4(n+m-1)]=(1+\rho)(1/4)\;,$
  \item[] $d_2=(3n+m-2)/[4(n+m-1)]=(2r+1-\rho)/4\;,$
\end{itemize}
the result can be explained as follows (see Figure~\ref{diag}):
  \begin{itemize}
    \item If $0\leq d\leq d_1$ then $(S)$ admits a solution $(\lambda,x)$
      such that: $\lambda=d/(1+\rho)\;$,
    \item If $d_1\leq d\leq d_2$ the $(S)$ admits a solution $(\lambda,x)$
      such that: $\lambda=1/4\;$,
    \item If $d_2< d\leq r$ or $r\leq d< d_2$ which cases
      correspond respectively to $r>1/2$ or $r<1/2$ then $(S)$ admits a solution $(\lambda,x)$
      such that: $\lambda=(r-d)/(2r-1+\rho)\;$,
    \item If $r\leq d\leq 1$ then $(S)$ admits a solution $(\lambda,x)$
      such that: $\lambda=0\;$.
  \end{itemize}
\end{rema}

\proof ~
\begin{itemize}
    \item If $\;0\leq d\leq d_1$\;, then
       a solution $(\lambda ,x)$ is given by:
       $$\lambda=\frac{d}{1+\rho}=\frac{n+m-1}{n+m}\;d,$$
       $$\begin{bmatrix}x_n\\x_{n+m}\\x_{1}\\x_{n+1}\end{bmatrix}
           =\begin{bmatrix}
             b_n\mpdiv \lambda^{n-1}\\
             \lambda^{n+1}\mpdiv a_n^2\mpdiv b_n\\
             e\\
             a_{n+m}\mpdiv a_n
        \end{bmatrix},$$
  which is a solution of:
  $$\left\{ \begin{array}{l}
        x_n=(b_n\mpdiv \lambda^{n-1}) x_1 , \\
        x_{n+m}=(b_m\mpdiv \lambda^{m-1}) x_{n+1}, \\
        x_1=(a_n\mpdiv \lambda)\sqrt{x_n x_{n+m}}, \\
        x_{n+1}=(a_{n+m}\mpdiv\lambda)\sqrt{x_n x_{n+m}}.
  \end{array} \right.$$
  Indeed:\\~ \\
  $\begin{array}{lll}
         \left[(\bar{a}_n\mpdiv\lambda^2)x_1x_{n+1}\mpdiv x_{n+m}\right]\mpdiv x_n
           & = & 1\mpdiv \lambda^4, \quad
           \text{ because } \;\;\bar{a}_n=1\mpdiv (a_n a_{n+m}),\\ \\
           & \geq & e, \;\; \text{ because } \;\;\lambda=\frac{n+m-1}{n+m}\;d\leq\frac{n+m-1}{n+m}\; 1/4\leq 1/4.
       \end{array}$\\~ \\~ \\
  $\begin{array}{lll}
         \left[(\bar{a}_{n+m}\mpdiv\lambda)x_1x_{n+1}\mpdiv x_n\right]\mpdiv x_{n+m}
           & = & 1\mpdiv \lambda^3, \quad
           \text{ because } \;\;\bar{a}_{n+m}=1\mpdiv (a_n a_{n+m}),\\ \\
           & \geq & e, \;\; \text{ because } \;\;\lambda\leq 1/4<1/3.
       \end{array}$\\~ \\~ \\
  $\begin{array}{lll}
         \left[(\bar{b}_n\mpdiv\lambda^{n-1})x_n\right]\mpdiv x_1
           & = & (1\mpdiv \lambda^2)^{n-1}, \quad
           \text{ because } \;\;\bar{b}_n b_n=1^{n-1},\\ \\
           & \geq & e, \;\; \text{ because } \;\;\lambda\leq 1/4<1/2.
       \end{array}$\\~ \\~ \\
  $\begin{array}{lll}
         \left[(\bar{b}_m\mpdiv\lambda^{m-1})x_{n+m}\right]\mpdiv x_{n+1}
           & = & 1^{m-1}\mpdiv(b_nb_ma_na_{n+m})\lambda^{n-m+2}, \quad
           \text{ because } \;\;\bar{b}_m=1^{m-1}\mpdiv b_m,\\ \\
           & = & (1\mpdiv \lambda^2)^{m-1}, \quad
           \text{ because }
           b_nb_ma_na_{n+m}=d^{n+m-1}=\lambda^{n+m},\\ \\
           & \geq & e, \;\; \text{ because } \;\;\lambda\leq 1/4<1/2.
       \end{array}$\\~ \\~ \\
    \item If $d_1\leq d\leq d_2$ then
      a solution $(\lambda ,x)$ is given by:
      $$\lambda=1/4$$
       $$\begin{bmatrix}x_n\\x_{n+m}\\x_{1}\\x_{n+1}\end{bmatrix}
        =\begin{bmatrix}
           \lambda^{m-3}\bar{a}_n\mpdiv b_m\\
           b_m a_{n+m}\mpdiv a_n\mpdiv \lambda^{m-1}\\
           e\\
           a_{n+m}\mpdiv a_n
        \end{bmatrix},$$
  which is a solution of:
  $$\left\{ \begin{array}{l}
        x_n=(\bar{a}_n\mpdiv\lambda^2) x_1 x_{n+1}\mpdiv x_{n+m}, \\
        x_{n+m}=(b_m\mpdiv \lambda^{m-1}) x_{n+1}, \\
        x_1=(a_n\mpdiv \lambda)\sqrt{x_n x_{n+m}}, \\
        x_{n+1}=(a_{n+m}\mpdiv\lambda)\sqrt{x_n x_{n+m}}.
  \end{array} \right.$$
  Indeed:\\~ \\
  $\begin{array}{lll}
         \left[(b_n\mpdiv\lambda^{n-1})x_1\right]\mpdiv x_n
           & = & (a_na_{n+m}b_nb_m)\mpdiv 1\mpdiv \lambda^{n+m-4}, \quad
           \text{ because } \;\;\bar{a}_n=1\mpdiv (a_n a_{n+m}),\\ \\
           & = & d^{n+m-1}\mpdiv 1\mpdiv (1/4)^{n+m-4}, \\ \\
           &   & \text{ because }\;\;a_na_{n+m}b_nb_m=d^{n+m-1},\quad\text{and}\quad \lambda=1/4,\\ \\
           & \geq & \frac{n+m}{4}-1-\frac{n+m-4}{4}=e, \quad
           \text{ because } \;\;d\geq d_1\;.
       \end{array}$\\~ \\~ \\
  $\begin{array}{lll}
         \left[(\bar{a}_{n+m}\mpdiv \lambda)x_1x_{n+1}\mpdiv x_n\right]\mpdiv x_{n+m}
           & = & \lambda = 1/4 \geq e.
       \end{array}$\\~ \\~ \\
  $\begin{array}{lll}
         \left[(\bar{b}_n\mpdiv\lambda^{n-1})x_n\right]\mpdiv x_1
           & = & 1^n\mpdiv(a_na_{n+m}b_nb_m)\mpdiv \lambda^{m-n-2},\\ \\
           &   & \text{ because } \;\;\bar{b}_n=1^{n-1}\mpdiv b_n,
           \quad\text{and}\quad \bar{a}_n=1\mpdiv (a_n a_{n+m}),\\ \\
           & = & 1^n\mpdiv d^{n+m-1}\mpdiv \lambda^{n-m+2}, \quad
           \text{ because }\;\;a_na_{n+m}b_nb_m=d^{n+m-1},\\ \\
           & \geq & n-\frac{3n+m-2}{4}-\frac{n-m+2}{4} = e, \quad
           \text{ because }\;\;d\leq d_2.
       \end{array}$\\

  $\begin{array}{lll}
         \left[(\bar{b}_m\mpdiv\lambda^{m-1})x_{n+m}\right]\mpdiv x_{n+1}
           & = & (1\mpdiv\lambda^2)^{m-1}, \quad
           \text{ because } \;\;\bar{b}_m b_m=1^{m-1},\\ \\
           & = & \frac{m-1}{2}>e, \quad \text{ because } \;\;
           \lambda=1/4.
       \end{array}$
  \item If $d_2< d\leq r$ or $r\leq d< d_2$, then
      \begin{equation}\label{lam3}
        \lambda=\frac{r-d}{2r-1+\rho}=\frac{n}{n-m+2}-\frac{n+m-1}{n-m+2}\;d
      \end{equation}
       $$\begin{bmatrix}x_n\\x_{n+m}\\x_{1}\\x_{n+1}\end{bmatrix}
        =\begin{bmatrix}
           \lambda^{n-1}\mpdiv \bar{b}_n\\
           b_m^2 a_{n+m}^2\mpdiv \bar{b}_n\mpdiv \lambda^{2m-n+1}\\
           e\\
           b_m a_{n+m}^2\mpdiv \bar{b}_n\mpdiv \lambda^{m-n+2}
        \end{bmatrix},$$
    which is a solution of:
    $$\left\{ \begin{array}{l}
      x_n=(\bar{a}_n\mpdiv\lambda^2) x_1 x_{n+1}\mpdiv x_{n+m}\;,\\
      x_{n+m}=(b_m\mpdiv\lambda^{m-1}) x_{n+1}\;,\\
      x_1=(\bar{b}_n\mpdiv\lambda^{n-1}) x_n\;,\\
      x_{n+1}=(a_{n+m}\mpdiv\lambda) \sqrt{x_n x_{n+m}}\;,
    \end{array} \right.$$
Indeed, $\lambda$ given by (\ref{lam3}) satisfies $0\leq \lambda\leq
1/4$ because:
     \begin{enumerate}
       \item If $d_2\leq d\leq r$, which corresponds to $n-m+2>0$, then
         $\frac{n+m-1}{n-m+2}>0$, and we can check that:
         $d_2\leq d\leq r\;\Rightarrow\;0\leq \lambda \leq 1/4$.
       \item If $r\leq d \leq d_2$, which corresponds to $n-m+2<0$, then
         $\frac{n+m-1}{n-m+2}<0$, and we can check that:
         $r\leq d\leq d_2\;\Rightarrow\;0\leq \lambda \leq 1/4$.
     \end{enumerate}
Then:\\~\\

  $\begin{array}{lll}
         \left[(b_n\mpdiv \lambda^{n-1})x_1\right]\mpdiv x_n
           & = & \left[1\mpdiv \lambda^2\right]^{n-1}, \quad
                 \text{ because } \;\;b_n \bar{b}_n=1^{n-1},\\ \\
           & \geq & e,\quad \text{ because } \;\;\lambda\leq 1/4<1/2.
       \end{array}$\\~

  $\begin{array}{lll}
         \left[(\bar{a}_{n+m}\mpdiv\lambda)x_1x_{n+1}\mpdiv x_n\right]\mpdiv x_{n+m}
           & = & 1^n\lambda^{m-n-1}\mpdiv a_n\mpdiv a_{n+1}\mpdiv b_m\mpdiv b_n,\\ \\
           &   & \text{ because } \;\;\bar{a}_{n+m}=1\mpdiv a_n\mpdiv a_{n+m}
             \quad \text{and} \quad \bar{b}_n=1^{n-1}\mpdiv b_n,\\ \\
           & = & (1^n\lambda^{n-m+2}\mpdiv d^{n+m-1})\lambda\;,\\ \\
           &   & \text{ because } \;\;a_n a_{n+m}b_mb_n=d^{n+m-1}\;,\\ \\
           & = & \lambda\geq e\;, \quad \text{because} \quad
                   \lambda^{n-m+2}=1^n\mpdiv d^{n+m-1}.
       \end{array}$\\~

  $\begin{array}{lll}
         \left[(a_n\mpdiv \lambda)\sqrt{x_nx_{n+m}}\right]\mpdiv x_1
           & = & a_na_{n+m}b_nb_m\mpdiv 1^{n-1}\lambda^{n-m-2},\quad
                 \text{ because } \;\;\bar{b}_n=1^{n-1}\mpdiv b_n,\\ \\
           & = & (d^{n+m-1}\mpdiv 1^n \lambda^{n-m+2}) (1\mpdiv
                 \lambda^4)\;,\;
                 \text{ because } \;\;a_n a_{n+m}b_mb_n=d^{n+m-1}\;,\\ \\
           & = & 1\mpdiv \lambda^4 \quad \text{because} \quad
                   \lambda^{n-m+2}=1^n\mpdiv d^{n+m-1},\\ \\
           & \geq & e,\quad \text{ because }\;\;\lambda\leq 1/4.
       \end{array}$\\~

  $\begin{array}{lll}
         \left[(\bar{b}_m\mpdiv \lambda^{m-1})x_{n+m}\right]\mpdiv x_{n+1}
           & = & 1^{m-1}\mpdiv \lambda^{2m-2},\quad
                 \text{ because } \;\;b_m\bar{b}_m=1^{m-1}\;\;,\\ \\
           & = & (1\mpdiv\lambda^2)^{m-1},\\ \\
           & \geq & e,\quad \text{ because }\;\;\lambda\leq 1/4<1/2.
       \end{array}$
  \item If $r\leq d\leq 1$, then $\lambda=0$

       $$\begin{bmatrix}x_n\\x_{n+m}\\x_{1}\\x_{n+1}\end{bmatrix}
        =\begin{bmatrix}
           e\mpdiv \bar{b}_n\\
           1^{n+1}\mpdiv a_n^2\mpdiv b_n\\
           e\\
           1a_{n+m}\mpdiv a_n
        \end{bmatrix},$$
    which is a solution of:
    $$\left\{ \begin{array}{l}
      x_n=(\bar{a}_n\mpdiv\lambda^2) x_1 x_{n+1}\mpdiv x_{n+m}\;,\\
      x_{n+m}=(\bar{a}_{n+m}\mpdiv\lambda) x_1 x_{n+1}\mpdiv x_n\;,\\
      x_1=(\bar{b}_n\mpdiv\lambda^{n-1}) x_n\;,\\
      x_{n+1}=(a_{n+m}\mpdiv\lambda) \sqrt{x_n x_{n+m}}\;,
    \end{array} \right.$$
Indeed;\\~ \\
  $\begin{array}{lll}
         \left[(b_n\mpdiv\lambda^{n-1})x_1\right]\mpdiv x_n
           & = & 1^{n-1}\;\geq e\;,\quad
                 \text{ because } \;\;b_n\bar{b}_n=1^{n-1},
       \end{array}$\\~

  $\begin{array}{lll}
         \left[(b_m\mpdiv\lambda^{m-1})x_{n+1}\right]\mpdiv x_{n+m}
           & = & a_na_{n+m}b_nb_m\mpdiv 1^n\;,\\ \\
           & = & d^{n+m-1}\mpdiv 1^n\;, \quad \text{ because }\;\; a_na_{n+m}b_nb_m=d^{n+m-1},\\ \\
           & \geq & e,\quad \text{ because }\;\;d\geq r.
       \end{array}$\\~

  $\begin{array}{lll}
         \left[(a_n\mpdiv \lambda)\sqrt{x_nx_{n+m}}\;\right]\mpdiv x_1
           & = & 1\;\geq e\;,
       \end{array}$\\~

  $\begin{array}{lll}
         \left[(\bar{b}_m\mpdiv \lambda^{m-1})x_{n+m}\right]\mpdiv x_{n+1}
           & = & 1^{n+m-1}\mpdiv (a_na_{n+m}b_nb_m)\;,\quad
                 \text{ because } \;\;\bar{b}_m=1^{m-1}\mpdiv b_m,\\  \\
           & = & 1^{n+m-1}\mpdiv d^{n+m-1}\;,\quad
                 \text{ because } \;\; a_na_{n+m}b_nb_m=d^{n+m-1},\\ \\
           & \geq & e\;, \quad \text{ because } d\leq 1\text{\finprv}
       \end{array}$
\end{itemize}

\begin{coro}\label{fond-cor1}
  In the case $r\geq 1/2$, a non negative eigenvalue $\lambda$ of $(S)$ is
  given by:
  $$\lambda=\max\left\{\min\left\{\frac{1}{1+\rho}\;d\;,\frac{1}{4}\;,\frac{r-d}{2r-1+\rho}\right\},0\right\}.$$
\end{coro}
\proof follows directly from Theorem~\ref{eigen}\finprv

\begin{coro}\label{fond-cor2}
  For large values of $n$ and $m$ such that $n>m-2$ (which is the case $r\geq 1/2$),
  a non negative eigenvalue $\lambda$ of (S) is given by:
  $$\lambda=\max\left\{\min\left\{d,\;\frac{1}{4},\;\frac{r-d}{2r-1}\right\},0\right\}.$$
\end{coro}
\proof follows directly from Corollary~\ref{fond-cor2}\finprv

\begin{rema}\label{rem1}
  We can check that as soon as we assume $m>1$, we get $d_1<d_2$, so
  we have: $0<d_1<d_2<1$. The position of $r$ with respect to $d_1$
  and $d_2$ gives three cases and divides the interval $[0,1]$, in each case
  in four regions, where the eigenvalues $\lambda$ of $(S)$ satisfy (see Figure \ref{diag}):
  \begin{itemize}
    \item[A.] $d\in[0,\min(d_1,r)[ \quad \Rightarrow \quad \lambda=\frac{1}{1+\rho}\;d$,
    \item[B.] $d\in[\min(d_1,r),d_1[ \quad \Rightarrow \quad
      \begin{cases}
        & \lambda=\frac{1}{1+\rho}\;d,\\
        \text{or} &
        \lambda=\frac{r-d}{2r-1+\rho},\\
        \text{or} & \lambda=0.
      \end{cases}$,
    \item[C.] $d\in[d_1,\min(d_2,r)[ \quad \Rightarrow \quad \lambda=\frac{1}{4}$,
    \item[D.] $d\in[\max(d_1,r),d_2[ \quad \Rightarrow \quad
      \begin{cases}
        & \lambda=\frac{1}{4},\\
        \text{or} &
        \lambda=\frac{r-d}{2r-1+\rho},\\
        \text{or} & \lambda=0.
      \end{cases}$,
    \item[E.] $d\in[d_2,\max(d_2,r)[ \quad \Rightarrow \quad \lambda=\frac{r-d}{2r-1+\rho}$,
    \item[F.] $d\in[\max(d_2,r),1] \quad \Rightarrow \quad \lambda=0$.
  \end{itemize}
\end{rema}

%\begin{figure}[htbp]
%  \begin{center}
%    \includegraphics[width=14cm,height=3.5cm]{diagrams}
%    \caption{The curve of $\lambda$ (given in Theorem~\ref{eigen}) depending on $d$.}
%    \label{diag}
%  \end{center}
%\end{figure}

\begin{lemm}\label{th-pos}
  If $\lambda>0$ then the system $(S)$ (solved on $(\lambda,x)$) is
  equivalent to the following system (solved on $(\lambda,z)$):
      \begin{align}
        & z_n=(\bar{a}_n\mpdiv b_m\mpdiv \lambda^2) z_1 \oplus (b_n\mpdiv
             \lambda^{n+m-2}) z_1, \label{z2}\\
        (SZ):\quad \quad \quad
        & z_{n+m}=b_m z_{n+1}, \label{z4}\\
        & z_1=(a_n\mpdiv\lambda)\sqrt{z_n z_{n+m}}\oplus (\bar{b}_n\mpdiv
            \lambda^{n-m}) z_n. \label{z1}\\
        & z_{n+1}= (a_{n+m}\mpdiv \lambda)\sqrt{z_n z_{n+m}} \oplus (\bar{b}_m\mpdiv
            \lambda^{2m-2}) z_{n+m}.\label{z3}
      \end{align}
\end{lemm}

\proof From the equations (\ref{S-nimplicit}) and (\ref{S-implicit})
we obtain:
$$\begin{cases}
  x_n x_{n+m}\leq (\bar{a}_n\mpdiv \lambda^2) x_1 x_{n+1}, &\\
  x_n x_{n+m}\leq (\bar{a}_{n+m}\mpdiv\lambda) x_1 x_{n+1}.&
\end{cases}$$
Since $\lambda >0$ we have: $x_n x_{n+m}< (\bar{a}_{n+m}\mpdiv
\lambda)x_1 x_{n+1}$. Thus:
\begin{equation}\label{eq-01}
  x_{n+m}=(b_m\mpdiv\lambda^{m-1}) x_{n+1}.
\end{equation}
By replacing $x_{n+m}$ in (\ref{S-nimplicit}), we obtain:
\begin{equation}\label{eq-02}
  x_n=\big[\bar{a}_n\mpdiv b_m \lambda^{m-3} \oplus
       b_n\mpdiv\lambda^{n-1}\big] x_1.
\end{equation}
The system $(S)$ is then equivalent to the system \{(\ref{eq-02}),
(\ref{eq-01}), (\ref{S-ar1}), (\ref{S-ar2})\}. On this later system
we use the following changing variable:
\begin{align}
  & z_n=x_n,  \nonumber \\
  & z_{n+m}=x_{n+m}\lambda^{2m-2},, \nonumber \\
  & z_1=x_1 \lambda^{m-1}, \nonumber \\
  & z_{n+1}=x_{n+1}\lambda^{m-1}, \nonumber
\end{align}
and we obtain the system $(SZ)$\finprv

%Let us define, when it exists the average growth rate vector $\chi$
%of the dynamical system $(DS)$ by:
%$$\chi=\lim_{k\to +\infty} \frac{x^k}{k}\;.$$

\begin{theo}
   If $r>1/2$ (that is $n\geq m$), and for the densities $d$
   satisfying $0<d<r$, the system $(S)$, and thus the eigenvalue problem
   $(EV)$, admit a unique positive eigenvalue $\lambda$
%   and the dynamical system $(DS)$ admits a unique average growth rate
%   vector $\chi$. Moreover we have: $\chi_i=\lambda,\;\forall 1\leq i\leq
%   n+m$, and $\lambda$ is
    given by:
   \begin{equation}\label{lamb}
      \lambda=\min\big\{\frac{1}{1+\rho}\;d\;,\;\frac{1}{4}\;,\;\frac{r-d}{2r-1+\rho}\big\}\;>0\;.
   \end{equation}
   This situation corresponds to the phases A, C and E of the case 3
   on Figure~\ref{diag} of Remark \ref{rem1}.
\end{theo}

\proof~
\begin{itemize}
  \item Let $\lambda$ be positive. Lemma \ref{th-pos} gives the
    equivalence of the systems $(S)$ and $(SZ)$. The later system
    is the eigenvalue problem associated to the following dynamical system:
    \begin{align}
        & z_n^k=(\bar{a}_n\mpdiv b_m) z_1^{k-2} \oplus b_n
             z_1^{k-(n+m-2)}, \label{zz2}\\
        & z_{n+m}^k=b_m z_{n+1}^k, \label{zz4}\\
        & z_1^k=a_n\sqrt{z_n^{k-1} z_{n+m}^{k-1}}\oplus \bar{b}_n
            z_n^{k-(n-m)}. \label{zz1}\\
        & z_{n+1}^k= a_{n+m}\sqrt{z_n^{k-1} z_{n+m}^{k-1}} \oplus \bar{b}_m
            z_{n+m}^{k-(2m-2)}.\label{zz3}
      \end{align}
    If $r>1/2$, that is $n\geq m$, then this dynamical system is implicit but triangular.
    Indeed, an iteration of the dynamics is to compute $z_n^k$ and $z_{n+1}^k$ in
    parallel, then compute $z_{n+m}^k$, and finally compute $z_1^k$.
    So the system $(SZ)$ can be
    interpreted as a dynamic programming equation of a stochastic
    optimal control problem, where $\lambda$ is the average optimal cost
    by unit of time. Since $\lambda$ is supposed to be positive, and
    from Corollary \ref{fond-cor1}, we obtain (\ref{lamb}).
  \item Let $\lambda_1>0$ and $\lambda_2>0$ be two positive
    eigenvalues of $(S)$. Lemma \ref{th-pos} tells us that both of
    $\lambda_1$ and $\lambda_2$ are eigenvalues of the system $(SZ)$.
    Since in the case when $r>2$, the system $(SZ)$ is a dynamic
    programming equation of a stochastic optimal control problem,
    which thus admits a unique eigenvalue, we conclude that
    $\lambda_1=\lambda_2$\finprv
\end{itemize}

%------------------------
\section{Conclusion}
%------------------------

The results of this note give a solution to the additive
eigenvalue problem associated to a dynamics of an elementary
2D-traffic system (two circular roads crossing at one junction,
managed by the priority-to-the-right-rule). The eigenvalue
$\lambda$, which is not necessarily unique, is given as a function
of two main quantities which are interpreted in terms of traffic
as the density $d$ of vehicles in the system, and the ratio $r$
between the non priority road size and the size of the whole
system. Moreover, we showed that when $r$ satisfies $r>1/2$, that
is when the size of the non priority road is bigger than the size
of the priority road, the uniqueness of $\lambda$, which is
positive in this case, is proved for densities satisfying $d<r$.


\begin{thebibliography}{50}

\bibitem{BCOQ92}
F.~Baccelli, G.~Cohen, G.J. Olsder, and J.P. Quadrat :
\emph{Synchronization and Linearity}, Wiley, 1992.

\bibitem{Far08} N. Farhi, \emph{Mod\'elisation Minplus et Commande du Trafic de
  Villes R\'eguli\`eres}, PhD Thesis Paris~1 University, 2008.

\bibitem{FGQ05} N.  Farhi, M. Goursat,  J.-P. Quadrat~: \emph{Derivation
of  the fundamental  traffic  diagram  for two  circular  roads
and a crossing using minplus algebra and Petri net modeling}, in
Proceedings of the 44th IEEE - CDC, Sevilla, 2005.

\bibitem{LMQ01} P. Lotito, E. Mancinelli and J.P. Quadrat \emph{A Minplus
Derivation of the Fundamental Car-Traffic Law},  IEEE Transactions
on Automatic Control V.50, N.5, p.699-705 May 2005.

\end{thebibliography}
\end{document}